\documentclass{amsart}
\usepackage{amsmath,amssymb}
\usepackage{yhmath}

\newcommand{\bdis}{\begin{displaymath}}
\newcommand{\edis}{\end{displaymath}}
\newcommand{\be}{\begin{equation}}
\newcommand{\ee}{\end{equation}}
\newcommand{\mbb}{\mathbb}
\newcommand{\mcal}{\mathcal}

\newcommand{\vp}{\varphi} 
\newcommand{\vth}{\vartheta}

\newcommand{\zf}{\zeta\left(\frac{1}{2}+it\right)}

\newcommand{\zfne}{\zeta\left(\frac{1}{2}+it_{2\nu}\right)} 
\newcommand{\zfno}{\zeta\left(\frac{1}{2}+it_{2\nu+1}\right)}

\newcommand{\FR}{\frac{x^n+y^n}{z^n}}

\DeclareMathOperator{\im}{Im}

\theoremstyle{definition}

\theoremstyle{remark}
\newtheorem{remark}[]{Remark}

\newtheorem*{mydef11}{{\bf Theorem 1}}

\newtheorem*{mydef12}{{\bf Theorem 2}}

\newtheorem*{mydef13}{{\bf Theorem 3}}

\newtheorem*{mydef41}{{\bf Corollary 1}}

\newtheorem*{mydef51}{{\bf Lemma 1}}

\newtheorem*{mydef52}{{\bf Lemma 2}}

\newtheorem*{mydef53}{{\bf Lemma 3}}

\newtheorem*{mydef81}{{\bf Property 1}}

\newtheorem*{mydef82}{{\bf Property 2}}

\newtheorem*{mydefR1}{{\bf Result 1}} 
\newtheorem*{mydefR2}{{\bf Result 2}} 
\newtheorem*{mydefR3}{{\bf Result 3}} 
\newtheorem*{mydefR4}{{\bf Result 4}}

\numberwithin{equation}{section}

\begin{document}

\title[Jacob's ladders, our asymptotic formulae (1981)  \dots]{Jacob's ladders, our asymptotic formulae (1981) and next $\zeta$-equivalents of the Fermat-Wiles theorem together with decomposition and synthesis of the Riemann's $\zeta$-oscillators}

\author{Jan Moser}

\address{Department of Mathematical Analysis and Numerical Mathematics, Comenius University, Mlynska Dolina M105, 842 48 Bratislava, SLOVAKIA}

\email{jan.mozer@fmph.uniba.sk}


\begin{abstract}
In this paper we obtain new  $\zeta$-equivalents of the Fermat-Wiles theorem. These are generated by our asymptotic formulae (1981) which brought $33.3\%$ improvement of the Hardy-Littlewood exponent $\frac 14$ dated 1918. 
\end{abstract}
\maketitle

\section{Introduction} 

\subsection{} 

In the classical book \cite{26} Titchmarsh obtained interesting asymptotic formulae for the segment $[T,T+H_1]$ with   
\be \label{1.1} 
H_1=T^{3/4}\psi(T)\sqrt{\ln T}, 
\ee 
and $\psi=\psi(T)$ is a function increasing arbitrarily slowly to $+\infty$. These formulae give the results: 
\be \label{1.2} 
\begin{split}
& \sum_{T\leq t_{2\nu}\leq T+H_1}Z(t_{2\nu})=\frac{1}{2\pi}H_1\ln\frac{T}{2\pi}+\mcal{O}(T^{3/4}\ln^{3/2}T), \\ 
& \sum_{T\leq t_{2\nu+1}\leq T+H_1}Z(t_{2\nu+1})=-\frac{1}{2\pi}H_1\ln\frac{T}{2\pi}+\mcal{O}(T^{3/4}\ln^{3/2}T). 
\end{split}
\ee 
In \cite{6} -- \cite{9} we have obtained the following result in this direction: 
\be \label{1.3} 
\begin{split}
& \sum_{T\leq t_{2\nu}\leq T+H_2}Z(t_{2\nu})=\frac{1}{2\pi}H_2\ln\frac{T}{2\pi}+\mcal{O}(T^{1/6}\psi\ln^{5}T), \\ 
& \sum_{T\leq t_{2\nu+1}\leq T+H_2}Z(t_{2\nu+1})=-\frac{1}{2\pi}H_2\ln\frac{T}{2\pi}+\mcal{O}(T^{1/6}\psi\ln^{5}T), 
\end{split}
\ee 
where 
\be \label{1.4} 
H_2=T^{1/6}\psi\ln^5T,\ 0<\psi(T)<\frac 16\frac{\ln T}{\ln\ln T}. 
\ee 

\begin{mydefR1}
The exponent $\frac 16$ in the formula (\ref{1.3}) is a $77.7\%$ improvement of the exponent $\frac 34$ in the local formula (\ref{1.2}) of Titchmarsh. Let us notice, that we have obtained this result 30 years after the classical book of Titchmarsh \cite{26} was published.  
\end{mydefR1} 

Next, it follows from (\ref{1.3}) that the interval 
\be \label{1.5} 
\left[\frac 12+iT,\ \frac 12+i(T+T^{1/6}\psi^2\ln^5T)\right],\ T\geq T_0(\psi), 
\ee  
contains a zero of odd order of the Riemann's function $\zeta(s)$. 

\begin{mydefR2}
The exponent $\frac 16$ in (\ref{1.5}) is a $33.3\%$ improvement of the Hardy-Littlewood exponent $\frac 14$ dated 1918, see \cite{1}. 
\end{mydefR2} 

Let us remind explicitly, that we have obtained this result 83 years after the classical memoir of Hardy and Littlewood. 

\begin{remark}
If we, for example, put in $H_2$ (see (\ref{1.4})) 
\be \label{1.6} 
\psi=\sqrt{\ln T}, 
\ee 
then we obtain 
\be \label{1.7} 
H_2=T^{1/6}\ln^6T, 
\ee  
and the interval 
\bdis 
\left[\frac 12+iT,\ \frac 12+i(T+T^{1/6}\ln^6T)\right],\ T\geq T_0(\sqrt{\ln T}). 
\edis 
\end{remark} 

\subsection{} 

Next, if $\frac 12+i\gamma$ denotes the zero of odd order of Riemann's zeta-function $\zeta(s),\ s=\sigma+it$ on the critical line $\sigma=\frac 12$, then our Result 1 gives the following estimate: 
\be \label{1.8} 
\gamma''-\gamma'<A(\gamma')^{1/6}(\ln\gamma')^{6}
\ee 
for the distance of two consecutive zeros 
\be \label{1.9} 
\frac 12+i\gamma',\ \frac 12+i\gamma'';\ \gamma'<\gamma'' 
\ee  
of odd order on the critical line. 

Of course, there are small improvements of the exponent $\frac 16$ in the formula (\ref{1.8}): 
\begin{itemize}
	\item[(a)] Karatsuba, see \cite{4}, 
	\be \label{1.10} 
	\frac{5}{32}=0.15625=6.25\% \ \mbox{improvement}. 
	\ee  
	\item[(b)] Ivi\` c, see \cite{3}, 
	\be \label{1.11} 
	0.1559458\dots \approx 0.19\% \ \mbox{improvement}
	\ee 
	and so on. 
\end{itemize} 
However, Littlewood, long before the estimates (\ref{1.8}), (\ref{1.10}) and (\ref{1.11}), has proved in 1924 that, by the Riemann hypothesis, it is true that\footnote{See \cite{5}.} 
\be \label{1.12} 
\gamma''-\gamma'<\frac{A}{\ln\ln\gamma'}, 
\ee  
and this estimate is independent on the order of the zeros (\ref{1.9}). 

\begin{remark}
Wee see, that mentioned estimate (\ref{1.8}), and of course the exponents (\ref{1.10}) and (\ref{1.11}), are extremely distant from the expected result. That is the estimates lap up the others but the formula continues. In this paper we obtain, after 44 years, the point of contact between our formulae (\ref{1.3}) and the Fermat-Wiles theorem without whatever improvements of the type (\ref{1.10}), (\ref{1.11}). 
\end{remark} 

\subsection{} 

We show in this paper that our asymptotic formula (\ref{1.3}) from 1981, connected with the Titchmarsh formula from 1951 and the Hardy-Littlewood estimate from 1918, also generates new functionals. For example the following one: 
\be \label{1.13} 
\begin{split}
& \lim_{\tau\to\infty}\frac{1}{\tau}
\left\{
\sum_{t_{2\nu}\geq (4\pi\alpha\tau)^6}^{t_{2\nu}\leq (4\pi\alpha\tau)^6+H((4\pi\alpha\tau)^6)}\zeta\left(\frac 12+i\tau_{2\nu}\right)\right\}\times \\ 
& \left\{
\int_{[(4\pi\alpha\tau)^{6}]^ 7}^{[(4\pi\alpha\tau)^6+2l]^7}\prod_{r=0}^6\left|\zeta\left(\frac{1}{2}+i\vp_1^r(t)\right)\right|^2{\rm d}t
\right\}^{-1}=\frac{\alpha}{l} 
\end{split}
\ee 
for every fixed $\alpha,l>0$, where 
\be \label{1.14} 
H(w)=w^{1/6}\ln w, 
\ee  
and, of course, 
\be \label{1.15} 
[G]^7=\vp_1^{-7}(G). 
\ee 
Let the symbol 
\be \label{1.16} 
\FR,\ x,y,z,n\in\mbb{N},\ n\geq 3 
\ee 
denote the Fermat's rationals. Then it follows from (\ref{1.13}) that the $\zeta$-condition 
\be \label{1.17} 
\begin{split}
& \lim_{\tau\to\infty}\frac{1}{\tau}
\left\{
\sum_{t_{2\nu}\geq (4\pi\alpha\tau)^6}^{t_{2\nu}\leq (4\pi\alpha\tau)^6+H((4\pi\alpha\tau)^6)}\zeta\left(\frac 12+i\tau_{2\nu}\right)\right\}\times \\ 
& \left\{
\int_{[(4\pi\alpha\tau)^{6}]^ 7}^{[(4\pi\alpha\tau)^6+2\FR]^7}\prod_{r=0}^6\left|\zeta\left(\frac{1}{2}+i\vp_1^r(t)\right)\right|^2{\rm d}t
\right\}^{-1}\not=\alpha
\end{split}
\ee 
on the set of all Fermat's rationals and for every fixed 
\be \label{1.18} 
\alpha>0,\ \alpha\notin \{\FR\}
\ee 
(for example every irrational number) expresses the next $\zeta$-equivalent of the Fermat-Wiles theorem. 

\subsection{} 

In this paper we give also next formula of a new type. Namely, we prove that the following formula holds true ($\tau\to+\infty$): 
\be \label{1.19} 
\begin{split}
& \sum_{t_{2\nu}\geq \tau^6}^{t_{2\nu}\leq\tau^6+H(\tau^6)}\zfne\sim \\ 
& \left\{
\int_\tau^{\overset{1}{\tau}}\left|\zf\right|^2{\rm d}t
\right\}\times 
\left\{
\int_{[\tau^6]^7}^{[\tau^6+(2\pi(1-c))^{-1}]^7}\prod_{r=0}^6\left|\zeta\left(\frac{1}{2}+i\vp_1^r(t)\right)\right|^2{\rm d}t
\right\}.
\end{split}
\ee 
This asymptotic formula \emph{controls}, in the limit $\tau\to+\infty$, $\zeta$-oscillations on the set of segments 
\be \label{1.20} 
[\tau,\overset{1}{\tau}],\ [\tau^6,[\tau^6+H(\tau^6)]],\ [[\tau^6]^7,[\tau^6+(2\pi(1-c))^{-1}]]
\ee 
lying on the critical line $\sigma=\frac 12$ and separated by huge distances.

\section{Jacob's ladders: notions and basic geometrical properties}  

\subsection{}

In this paper we use the following notions of our works \cite{5} -- \cite{9}: 
\begin{itemize}
\item[{\tt (a)}] Jacob's ladder $\vp_1(T)$, 
\item[{\tt (b)}] direct iterations of Jacob's ladders 
\bdis 
\begin{split}
	& \vp_1^0(t)=t,\ \vp_1^1(t)=\vp_1(t),\ \vp_1^2(t)=\vp_1(\vp_1(t)),\dots , \\ 
	& \vp_1^k(t)=\vp_1(\vp_1^{k-1}(t))
\end{split}
\edis 
for every fixed natural number $k$, 
\item[{\tt (c)}] reverse iterations of Jacob's ladders 
\be \label{2.1}  
\begin{split}
	& \vp_1^{-1}(T)=\overset{1}{T},\ \vp_1^{-2}(T)=\vp_1^{-1}(\overset{1}{T})=\overset{2}{T},\dots, \\ 
	& \vp_1^{-r}(T)=\vp_1^{-1}(\overset{r-1}{T})=\overset{r}{T},\ r=1,\dots,k, 
\end{split} 
\ee   
where, for example, 
\be \label{2.2} 
\vp_1(\overset{r}{T})=\overset{r-1}{T}
\ee  
for every fixed $k\in\mbb{N}$ and every sufficiently big $T>0$. We also use the properties of the reverse iterations listed below.  
\be \label{2.3}
\overset{r}{T}-\overset{r-1}{T}\sim(1-c)\pi(\overset{r}{T});\ \pi(\overset{r}{T})\sim\frac{\overset{r}{T}}{\ln \overset{r}{T}},\ r=1,\dots,k,\ T\to\infty,  
\ee 
\be \label{2.4} 
\overset{0}{T}=T<\overset{1}{T}(T)<\overset{2}{T}(T)<\dots<\overset{k}{T}(T), 
\ee 
and 
\be \label{2.5} 
T\sim \overset{1}{T}\sim \overset{2}{T}\sim \dots\sim \overset{k}{T},\ T\to\infty.   
\ee  
\end{itemize} 

\begin{remark}
	The asymptotic behaviour of the points 
	\bdis 
	\{T,\overset{1}{T},\dots,\overset{k}{T}\}
	\edis  
	is as follows: at $T\to\infty$ these points recede unboundedly each from other and all together are receding to infinity. Hence, the set of these points behaves at $T\to\infty$ as one-dimensional Friedmann-Hubble expanding Universe. 
\end{remark}  

\subsection{} 

Let us remind that we have proved\footnote{See \cite{18}, (3.4).} the existence of almost linear increments 
\be \label{2.6} 
\begin{split}
& \int_{\overset{r-1}{T}}^{\overset{r}{T}}\left|\zf\right|^2{\rm d}t\sim (1-c)\overset{r-1}{T}, \\ 
& r=1,\dots,k,\ T\to\infty,\ \overset{r}{T}=\overset{r}{T}(T)=\vp_1^{-r}(T)
\end{split} 
\ee 
for the Hardy-Littlewood integral (1918) 
\be \label{2.7} 
J(T)=\int_0^T\left|\zf\right|^2{\rm d}t. 
\ee  

For completeness, we give here some basic geometrical properties related to Jacob's ladders. These are generated by the sequence 
\be \label{2.8} 
T\to \left\{\overset{r}{T}(T)\right\}_{r=1}^k
\ee 
of reverse iterations of the Jacob's ladders for every sufficiently big $T>0$ and every fixed $k\in\mbb{N}$. 

\begin{mydef81}
The sequence (\ref{2.8}) defines a partition of the segment $[T,\overset{k}{T}]$ as follows 
\be \label{2.9} 
|[T,\overset{k}{T}]|=\sum_{r=1}^k|[\overset{r-1}{T},\overset{r}{T}]|
\ee 
on the asymptotically equidistant parts 
\be \label{2.10} 
\begin{split}
& \overset{r}{T}-\overset{r-1}{T}\sim \overset{r+1}{T}-\overset{r}{T}, \\ 
& r=1,\dots,k-1,\ T\to\infty. 
\end{split}
\ee 
\end{mydef81} 

\begin{mydef82}
Simultaneously with the Property 1, the sequence (\ref{2.8}) defines the partition of the integral 
\be \label{2.11} 
\int_T^{\overset{k}{T}}\left|\zf\right|^2{\rm d}t
\ee 
into the parts 
\be \label{2.12} 
\int_T^{\overset{k}{T}}\left|\zf\right|^2{\rm d}t=\sum_{r=1}^k\int_{\overset{r-1}{T}}^{\overset{r}{T}}\left|\zf\right|^2{\rm d}t, 
\ee 
that are asymptotically equal 
\be \label{2.13} 
\int_{\overset{r-1}{T}}^{\overset{r}{T}}\left|\zf\right|^2{\rm d}t\sim \int_{\overset{r}{T}}^{\overset{r+1}{T}}\left|\zf\right|^2{\rm d}t,\ T\to\infty. 
\ee 
\end{mydef82} 

It is clear, that (\ref{2.10}) follows from (\ref{2.3}) and (\ref{2.5}) since 
\be \label{2.14} 
\overset{r}{T}-\overset{r-1}{T}\sim (1-c)\frac{\overset{r}{T}}{\ln \overset{r}{T}}\sim (1-c)\frac{T}{\ln T},\ r=1,\dots,k, 
\ee  
while our eq. (\ref{2.13}) follows from (\ref{2.6}) and (\ref{2.5}).  

\section{New $\zeta$-functionals and corresponding $\zeta$-equivalents of the Fermat-Wiles theorem generated by the formula (\ref{1.3})} 

\subsection{} 

Let us remind real-valued $Z$-function defined by Riemann\footnote{See \cite{23}, (35), (44), (62); \cite{24}, p. 98.} as 
\be \label{3.1} 
Z(t)=e^{i\vth(t)}\zf, 
\ee  
where 
\be \label{3.2} 
\begin{split}
& \vth(t)=-\frac t2\ln 2\pi+\im\ln\Gamma\left(\frac 14+i\frac t2\right)= \\ 
& \frac t2\ln\frac{t}{2\pi}-\frac t2-\frac{\pi}{8}+\mcal{O}\left(\frac 1t\right). 
\end{split}
\ee 
Next, the Gram's sequence $\{t_\nu\}$ is defined by the condition 
\be \label{3.3} 
\vth(t_\nu)=\pi\nu,\ \nu=1,2,\dots
\ee 
Since\footnote{See (\ref{3.1}), (\ref{3.3}).} 
\be \label{3.4} 
Z(t_{2\nu})=\zfne,\ Z(t_{2\nu+1})=-\zfno, 
\ee 
we can rewrite the formulae (\ref{1.13}) in the form 
\be \label{3.5} 
\begin{split}
& \sum_{T\leq t_{2\nu}\leq T+H_2}\zfne=\frac{1}{2\pi}H_2\ln\frac{T}{2\pi}+\mcal{O}(T^{1/6}\psi\ln^5T), \\ 
& \sum_{T\leq t_{2\nu+1}\leq T+H_2}\zfno=\frac{1}{2\pi}H_2\ln\frac{T}{2\pi}+\mcal{O}(T^{1/6}\psi\ln^5T). 
\end{split}
\ee  
It is sufficient to use the first formula (\ref{3.5}), that is\footnote{See (\ref{1.4}), (\ref{1.6}), (\ref{1.7}); $H_2=H$.} 
\be \label{3.6} 
\begin{split}
& \sum_{T\leq t_{2\nu}\leq T+H}\zfne=\frac{1}{2\pi}H\ln\frac{T}{2\pi}+\mcal{O}(T^{1/6}\ln^{11/2}T), \\ 
& H=H(T)=T^{1/6}\ln T, 
\end{split}
\ee 
or in the form 
\be \label{3.7} 
\sum_{T\leq t_{2\nu}\leq T+H(T)}\zfne=\frac{1}{2\pi}T^{1/6}\ln^7T+\mcal{O}(T^{1/6}\ln^6T), 
\ee 
and, finally, we get 

\begin{mydef51}
\be \label{3.8} 
\begin{split}
& \sum_{t_{2\nu}\geq T}^{t_{2\nu}\leq T+H(T)}\zfne= \\ 
& \frac{1}{2\pi}T^{1/6}\ln^7T\left\{1+\mcal{O}\left(\frac{1}{\ln T}\right)\right\},\ H(T)=T^{1/6}\ln T. 
\end{split}
\ee 
\end{mydef51} 

\subsection{} 

Next, we use in (\ref{3.8}) our formula\footnote{See \cite{21}, (3.18), $f_m=1$.} 
\be \label{3.9} 
\begin{split}
& \int_{\overset{k}{T}}^{\overset{k}{\wideparen{T+2l}}}\prod_{r=0}^{k-1}\left|\zeta\left(\frac 12+i\vp_1^r(t)\right)\right|^2{\rm d}t= \\ 
& \left\{
1+\mcal{O}_{k^2}\left(\frac{\ln\ln T}{\ln T}\right)
\right\}2l\ln^kT 
\end{split}
\ee 
for every fixed $l>0$ and $k\in \mbb{N}$ in the case $k=7$, i.e. the following formula 
\be \label{3.10} 
\begin{split}
	& 2l\ln^7T=  
	\left\{
	\int_{\overset{7}{T}}^{\overset{7}{\wideparen{T+2l}}}\prod_{r=0}^{6}\left|\zeta\left(\frac 12+i\vp_1^r(t)\right)\right|^2{\rm d}t
	\right\}\times 
	\left\{
	1+\mcal{O}\left(\frac{\ln\ln T}{\ln T}\right)
	\right\}. 
\end{split}
\ee  
Now we have the following result\footnote{See (\ref{3.8}), (\ref{3.10}).}.

\begin{mydef52}
\be \label{3.11} 
\begin{split}
& \left\{
\sum_{t_{2\nu}\geq T}^{t_{2\nu}\leq T+H(T)}\zfne
\right\}\times 
\left\{
\int_{\overset{7}{T}}^{\overset{7}{\wideparen{T+2l}}}\prod_{r=0}^{6}\left|\zeta\left(\frac 12+i\vp_1^r(t)\right)\right|^2{\rm d}t
\right\}^{-1} = \\ 
& \frac{1}{4\pi l}T^{1/6}\left\{
1+\mcal{O}\left(\frac{\ln\ln T}{\ln T}\right)
\right\} 
\end{split}
\ee 
for every fixed $l>0$. 
\end{mydef52} 

Consequently, the substitution 
\be \label{3.12} 
T=(4\pi\alpha\tau)^6,\ \{T\to+\infty\} \ \Leftrightarrow \ \{\tau\to+\infty\}
\ee 
(there are also other possibilities) in (\ref{3.11}) gives the following functional. 

\begin{mydef11}
\be\label{3.13} 
\begin{split}
& \lim_{\tau\to\infty}\frac{1}{\tau}
\left\{
\sum_{t_{2\nu}\geq (4\pi\alpha\tau)^6}^{t_{2\nu}\leq (4\pi\alpha\tau)^6+H((4\pi\alpha\tau)^6)}\zfne
\right\}\times \\ 
& \left\{
\int_{[(4\pi\alpha\tau)^6]^7}^{[(4\pi\alpha\tau)^6+2l]^7}\prod_{r=0}^6\left|\zeta\left(\frac 12+i\vp_1^r(t)\right)\right|^2{\rm d}t
\right\}^{-1}=\frac{\alpha}{l}
\end{split}
\ee 
for every fixed $\alpha,l>0$, where 
\bdis 
[Y]^7=\vp_1^{-7}(Y). 
\edis 
\end{mydef11}

In the special case of Fermat's rationals\footnote{See (\ref{1.16}).} 
\be \label{3.14} 
l\to\FR
\ee  
we obtain from (\ref{3.13}) the following result. 

\begin{mydef41}
\be \label{3.15} 
\begin{split}
& \lim_{\tau\to\infty}\frac{1}{\tau}
\left\{
\sum_{t_{2\nu}\geq (4\pi\alpha\tau)^6}^{t_{2\nu}\leq (4\pi\alpha\tau)^6+H((4\pi\alpha\tau)^6)}\zfne
\right\}\times \\
& \left\{
\int_{[(4\pi\alpha\tau)^6]^7}^{[(4\pi\alpha\tau)^6+2\FR]^7}\prod_{r=0}^6\left|\zeta\left(\frac 12+i\vp_1^r(t)\right)\right|^2{\rm d}t
\right\}^{-1}=\alpha\left(\FR\right)^{-1}
\end{split}
\ee 
for every fixed $l>0$ and every Fermat's rational. 
\end{mydef41} 

Consequently, we obtain from (\ref{3.15}) the next result. 
\begin{mydef12}
The $\zeta$-condition 
\be \label{3.16} 
\begin{split}
& \lim_{\tau\to\infty}\frac{1}{\tau}
\left\{
\sum_{t_{2\nu}\geq (4\pi\alpha\tau)^6}^{t_{2\nu}\leq (4\pi\alpha\tau)^6+H((4\pi\alpha\tau)^6)}\zfne
\right\}\times \\ 
& \left\{
\int_{[(4\pi\alpha\tau)^6]^7}^{[(4\pi\alpha\tau)^6+2\FR]^7}\prod_{r=0}^6\left|\zeta\left(\frac 12+i\vp_1^r(t)\right)\right|^2{\rm d}t
\right\}^{-1}\not=\alpha
\end{split}
\ee 
on the set of all Fermat's rationals and for every fixed 
\be \label{3.17} 
\alpha>0,\ \alpha\notin\left\{\FR\right\}
\ee 
expresses the next $\zeta$-equivalent of the Fermat-Wiles theorem. 
\end{mydef12} 

\begin{remark}
If we use the substitution 
\be \label{3.18}  
\alpha\to\FR
\ee 
instead of (\ref{3.14}), then we obtain next $\zeta$-equivalent instead of (\ref{3.16}). 
\end{remark}

\section{On decomposition and synthesis of some sets of the Riemann's oscillators} 

\subsection{} 

The Riemann-Siegel formula\footnote{See \cite{23}, p. 60.} 
\be \label{4.1} 
Z(t)=\sum_{n\leq\sqrt{\frac{t}{2\pi}}}\frac{2}{\sqrt{n}}\cos\{\vth(t)-t\ln n\}+\mcal{O}(t^{-1/4})
\ee 
belongs to Riemann and it was restored by C.L. Siegel from Riemann's manuscripts. 

Let us remind the local structure of this Riemann-Siegel formula, i.e. the spectral form of this one\footnote{See \cite{17}, (3.1) -- (3.8).} 
\be \label{4.2} 
\begin{split}
& Z(t)=\sum_{n\leq\tau(x)}\frac{2}{\sqrt{n}}\cos\{t\omega_n(x)+\psi(x)\}+R(x), \\ 
& \tau(x)=\sqrt{\frac{x}{2\pi}} , \\ 
& R(x)=\mcal{O}(x^{-1/4}), \\ 
& t\in [x,x+v],\ v=\in[0,x^{1/4}], 
\end{split}
\ee 
where the functions 
\be \label{4.3} 
\frac{2}{\sqrt{n}}\cos\{t\omega_n(x)+\psi(x)\}
\ee  
are Riemann's oscillators with: 
\begin{itemize}
	\item[(a)] Amplitude 
	\bdis 
	\frac{2}{\sqrt{n}}, 
	\edis 
	\item[(b)] incoherent phase constant 
	\bdis 
	\psi(x)=-\frac{x}{2}-\frac{\pi}{8}, 
	\edis 
	\item[(c)] non-synchronised local time 
	\bdis 
	t\in [x,x+v], 
	\edis 
	\item[(d)] local spectrum of cyclic frequencies 
	\bdis 
	\{\omega_n(x)\}_{n\leq\tau(x)},\ \omega_n(x)=\ln\frac{\tau(x)}{n}. 
	\edis 
\end{itemize}

\subsection{} 

Now, we go back to our Lemma 2. Namely, the substitution 
\be \label{4.4} 
T=\tau^6 
\ee  
gives the following formula 
\be \label{4.5} 
\begin{split}
& \sum_{t_{2\nu}\geq \tau^6}^{t_{2\nu}\leq \tau^6+H(\tau^6)}\zfne = \\ 
& \frac{1}{4\pi l}\left\{1+\mcal{O}\left(\frac{\ln\ln\tau}{\ln\tau}\right)\right\}\tau\int_{[\tau^6]^7}^{[\tau^6+2l]^7}\prod_{r=0}^6\left|\zeta\left(\frac 12+i\vp_1^r(t)\right)\right|^2{\rm d}t
\end{split}
\ee 
for every fixed $l>0$. 

Next, we use our almost linear formula\footnote{See \cite{18}, (3.4), (3.6) with $r=1$.}
\bdis 
\begin{split}
& \int_\tau^{\overset{1}{\tau}}\left|\zf\right|^2{\rm d}t=(1-c)\tau+\mcal{O}(\tau^{1/3+\delta})= \\ 
& (1-c)\tau\left\{1+\mcal{O}(\tau^{-2/3+\delta})\right\}, 
\end{split}
\edis 
i.e. 
\be \label{4.6} 
(1-c)\tau=\left\{1+\mcal{O}(\tau^{-2/3+\delta})\right\}\int_\tau^{\overset{1}{\tau}}\left|\zf\right|^2{\rm d}t , 
\ee 
where $\delta>0$ is sufficiently small. And, as a consequence, we obtain the next Lemma. 

\begin{mydef53}
\be \label{4.7} 
\begin{split}
& \sum_{t_{2\nu}\geq \tau^6}^{t_{2\nu}\leq \tau^6+H(\tau^6)}\zfne = \\ 
& \frac{1}{4\pi(1-c)l}\left\{1+\mcal{O}\left(\frac{\ln\ln\tau}{\ln\tau}\right)\right\}\int_\tau^{\overset{1}{\tau}}\left|\zf\right|^2{\rm d}t\times \\ 
& \int_{[\tau^6]^7}^{[\tau^6+2l]^7}\prod_{r=0}^6\left|\zeta\left(\frac 12+i\vp_1^r(t)\right)\right|^2{\rm d}t
\end{split}
\ee 
for every fixed $l>0$. 
\end{mydef53} 

Now, in the case 
\be \label{4.8} 
2l=\frac{1}{2\pi(1-c)}, 
\ee  
(the condition of an asymptotic equilibrium), we obtain the following Theorem. 

\begin{mydef13}
\be \label{4.9} 
\begin{split}
& \sum_{t_{2\nu}\geq \tau^6}^{t_{2\nu}\leq \tau^6+H(\tau^6)}\zfne\sim \\
& \left\{\int_\tau^{\overset{1}{\tau}}\left|\zf\right|^2{\rm d}t\right\}\times 
\left\{\int_{[\tau^6]^7}^{[\tau^6+2l]^7}\prod_{r=0}^6\left|\zeta\left(\frac 12+i\vp_1^r(t)\right)\right|^2{\rm d}t\right\}
\end{split}
\ee 
for $\tau\to\infty$. 
\end{mydef13} 

\begin{remark}
The nonlinear formula (\ref{4.9}) expresses \emph{the interactions} of three sets of Riemann's oscillators\footnote{See (\ref{3.4}), (\ref{4.2}) and (\ref{4.3}).} contained in the sum 
\be \label{4.10} 
\sum_{t_{2\nu}\geq \tau^6}^{t_{2\nu}\leq \tau^6+H(\tau^6)}Z(t_{2\nu})
\ee 
and also in two integrals 
\be \label{4.11} 
\int_\tau^{\overset{1}{\tau}}Z^2(t){\rm d}t,\ \int_{[\tau^6]^7}^{[\tau^6+2l]^7}\prod_{r=0}^6Z^2(\vp_1^r(t)){\rm d}t 
\ee 
(more complicated Riemann's oscillators are contained in the last two integrals). The nonlinear interactions, that proceed in the left-right direction, represent the decomposition (from simpler to more complicated) and that one in the opposite direction the synthesis (from complicated to simple structure). 
\end{remark} 

\section{The evolution of Titchmarsh discrete method (1934)} 

We have proved in our paper \cite{20} a new $\zeta$-equivalent of the Fermat-Wiles theorem. This proof is based on our asymptotic formulas connected with our proof (1980) of the Titchmarsh hypothesis (1934). 

In the present paper, we have proved next $\zeta$-equivalents. These are also based on our old formulas (\ref{1.3}) connected with the classical Titchmarsh formulas (1951). In connection with notices above, let us give more detailed remarks about the period 1976 -- 1981. 

\subsection{} 

First, we have added to our Result 1 and Result 2 also the following two. 

Let $N_0(T)$ denote the number of of zeros of the function 
\bdis 
\zf,\ t\in (0,T). 
\edis 
In 1921 Hardy and Littlewood have obtained the estimate \cite{2}: 
\be \label{5.1} 
N_0(T+T^{1/2+\epsilon})-N_0(T)> A(\epsilon)T^{1/2+\epsilon},\ T>T_0(\epsilon). 
\ee  
Later, in 1942, Selberg has derived the estimate \cite{22}: 
\be \label{5.2} 
N_0(T+T^{1/2+\epsilon})-N_0(T)>A(\epsilon)T^{1/2+\epsilon}\ln T,\ T\geq T_0(\epsilon), 
\ee 
and has conjectured that the exponent $1/2$ in (\ref{5.2}) can be replaced by a number less than $1/2$. 

After 62 years we have obtained, in the direction of Hardy-Littlewood result (\ref{5.1}), the following estimate\footnote{See \cite{11}, \cite{12}.}: 
\be \label{5.3} 
N_0(T+T^{5/12}\psi\ln^3T)-N_0(T)>A(\psi)T^{5/12}\psi\ln^3T,\ T\geq T_0(\psi)
\ee 
as a result of the synthesis of the following deep concepts due to English mathematicians: 
\begin{itemize}
	\item[(a)] The averaging method of Hardy and Littlewood. 
	\item[(b)] The discrete method of Titchmarsh, see \cite{24}. 
	\item[(c)] The Titchmarsh variant of the van der Corput's method for estimation of corresponding trigonometric sums, see \cite{25}. 
\end{itemize} 

\begin{mydefR3}
The exponent $5/12$ represents a $16.6\%$-improvement of the Hardy-Littlewood exponent $1/2$ and, simultaneously, the first step toward the proof of the Selberg's conjecture. 
\end{mydefR3} 

Next, let us remind that in 1934 Titchmarsh presented the following hypothesis\footnote{See \cite{24}, p. 105.}: There is such positive $A$ that 
\be \label{5.4} 
\sum_{\nu=M+1}^N Z^2(t_\nu)Z^2(t_{\nu+1})=\mcal{O}(N\ln^AN), 
\ee  
where (if $t_N=T$) 
\be \label{5.5} 
N\sim\frac{1}{2\pi}T\ln T, 
\ee  
and $M$ is a sufficiently big fixed number. 

\begin{mydefR4}
After 46 years, we have proved, see \cite{10}, the Titchmarsh hypothesis with $A=4$, i.e. 
\be \label{5.6} 
\sum_{\nu=M+1}^N Z^2(t_\nu)Z^2(t_{\nu+1})=\mcal{O}(N\ln^4N), \ N\to\infty. 
\ee 
\end{mydefR4} 

\begin{remark}
After next 45 years we have been able to prove the existence of next $\zeta$-equivalents  of the Fermat-Wiles theorem based on sums of Titchmarsh type, see \cite{20}. 
\end{remark} 

\begin{remark}
It is astonishing that in the whole period 1859 -- 2023 there was no indication of the fact that Fermat's last theorem, now Fermat-Wiles theorem, could have been connected with the properties of the Riemann zeta-function on the critical line $\sigma=\frac 12$. This was shown in our paper \cite{19}. 
\end{remark} 

\subsection{} 

On the basis of essential improvements\footnote{Comp. Results 1 -- 4.} of the classical results of Hardy and Littlewood (1918) and (1921), and also of the proof of Titchmarsh hypothesis (1934) I have been invited to the international conference in Moscow (Sept. 14 -- 19, 1981) organized in honour of the 90th. birthday of I.M. Vinogradov. 

\begin{remark}
This conference joined my research with that one of A.A. Karatsuba, namely: 
\begin{itemize}
	\item[(a)] I have presented the list of main results obtained in the period 1976 -- 1981, see \cite{13}\footnote{Comp. Results 1 -- 4 in this paper.}. 
	\item[(b)] A. A. Karatsuba published, see \cite{4}, the first work from his series of papers in the direction of my results, I sent him the separates as well as complete manuscripts of the papers accepted for publication. 
\end{itemize}
\end{remark} 

\begin{remark}
Our scientific positions were as follows: 
\begin{itemize}
	\item[(c)] A. A. Karatsuba was a full Professor and Member of the Steklov institute in Moscow. 
	\item[(d)] I was a lector without scientific degree (1966 -- 1990), studying individually the works of Hardy and Littlewood, Titchmarsh, Ingham and, of course, of Selberg, on the Riemann's zeta-function. 
\end{itemize}
\end{remark} 

\begin{remark}
Only after the fall of communism in Czechoslovakia in 1989 I have obtained directly the highest scientific degree D.Sc. on the basis of my dissertation \emph{The evolution of the Titchmarsh discrete method}, pp. I -- XIV, 1 -- 236, in Russian, unpublished. I am grateful for this to the kind initiative of prof. Stefan Schwarz, D. Sc., co-founder of the theory of semigroups and Member of the editorial board of the journal \emph{Semigroups forum}. 
\end{remark}

I would like to thank Michal Demetrian for his moral support of my study of Jacob's ladders.

\end{document}